\def\am            {\alpha^-} 
\def\ap            {\alpha^+} 
\def\apm           {\alpha^\pm} 
\def\be            {\begin{equation}}
\def\bearl         {\begin{array}{l}}
\def\bearll        {\begin{array}{ll}}
\def\bpm           {\beta^{+-}} 
\def\C             {\ensuremath{\mathcal C}}
\def\CA            {\ensuremath{\mathcal C_{\!A}}}
\def\CAA           {\ensuremath{\mathcal C_{\!A|A}}}
\def\CAB           {\ensuremath{\mathcal C_{\!A|B}}}
\def\cft           {conformal field theory}
\def\cfts          {conformal field theories}
\def\cir           {\,{\circ}\,}
\def\CM            {\ensuremath{\mathcal M}}
\def\Cs            {\ensuremath{\mathcal C_{\!\mathcal M}^*}}
\def\complex       {{\ensuremath{\mathbbm C}}}
\def\Cong          {\,{\cong}\,} 
\newcommand\D[3]   {D_{#1}^{#2#3}}  
\newcommand\dd[5]  {\mathrm d_{#1}^{~#2#3;#4\Ol{#5}}}

\newcommand\ddi[5] {\big(\mathrm{d}^{-1}_{}\big)^{~~~~#1}_{#2#3;#4\Ol{#5}}}
\def\dim           {\mathrm{dim}}
\def\dimc          {\dim_\complex}
\def\dsty          {\displaystyle }
\newcommand\e[3]   {e^{}_{#1;#2#3}}
\def\ee            {\end{equation}}

\def\eear          {\end{array}}
\newcommand\ei[4]  {e^{}_{(#1,#2);#3#4}}
\def\End           {\mathrm{End}}
\def\EndAA         {\mathrm{End}_{\!A|A}}
\def\eps           {\varepsilon}
\def\Eps           {\epsilon}
\def\eq            {\,{=}\,}
\newcommand\eqpic[4]{\begin{eqnarray}
                   \begin{picture}(#2,#3){}\end{picture}\nonumber\\
                   \raisebox{-#3pt}{ \begin{picture}(#2,#3) #4 \end{picture} }
                   \label{#1} \\~\nonumber \end{eqnarray} }
\newcommand\erf[1] {(\ref{#1})}
\def\Es            {\ensuremath{\mathcal E}}
\def\findim        {fi\-ni\-te-di\-men\-si\-o\-nal}
\def\Fs            {\ensuremath{\mathcal F}}
\newcommand\Ga[2]  {\Gamma_{\!#1#2}}
\newcommand\Gab[2] {\ol\Gamma_{\!#1#2}}
\def\Hom           {\mathrm{Hom}}
\def\HomA          {\mathrm{Hom}_{\!A}}
\def\HomAA         {\mathrm{Hom}_{\!A|A}}
\def\HomCs         {\mathrm{Hom}_{\Cs}}
\def\id            {\mbox{\sl id}}
\def\ida           {\id_{A}}
\def\idav          {\id_{A^\vee_{\phantom:}}}
\def\idM           {\id_{\dot M}}

\def\iN            {\,{\in}\,} 
\newcommand\Includeourbeautifulpicture[1] {{\begin{picture}(0,0)(0,0)
                  \scalebox{.38}{\includegraphics{pic_navf_#1.eps}}\end{picture}}}
\def\Is            {K}
\def\ka            {\kappa}
\def\kap           {{\kappa'}}
\def\kapp          {{\kappa''}}
\def\kc            {\ensuremath{K_0(\C)}}
\def\kcaa          {\ensuremath{K_0(\CAA)}}
\def\kcs           {\ensuremath{K_0(\Cs)}}
\newcommand\labl[1]{\label{#1}\ee}
\def\Mat           {\mathrm{Mat}}
\def\nE            {\,{\not=}\,}
\newcommand\Ns[3]  {N^{~~#3}_{#1#2}}
\def\obj           {{\mathcal O}bj}
\def\ol            {\overline }
\def\Ol            {} 
\def\one           {{\ensuremath{\mathbf 1}}}

\def\ota           {\,{\otimes_{\!A}}\,}
\def\Ota           {{\otimes_{\!A}}}
\def\oti           {\,{\otimes}\,}
\def\Oti           {{\otimes}}
\def\otim          {\,{\otimes^-}\,}
\def\otip          {\,{\otimes^+}} 
\def\oxi           {\,{\boxtimes}\,}
\def\Oxi           {{\boxtimes}}
\def\qed           {\hfill$\Box$}
\def\rep           {representation}
\def\rhol          {\rho_{\mathrm l}}
\def\rhor          {\rho_{\mathrm r}}
\def\rmd           {{\rm d}}
\def\sse           {\scriptsize}
\def\ssf           {nondegenerate}

\def\ssfy          {nondegeneracy}
\def\Ssfy          {Nondegeneracy}

\def\Star          {\,{\ast}\,}
\def\Times         {\,{\times}\,}
\def\To            {\,{\to}\,} 
\newcommand\ui[4]  {(u^{-1})^{#1}_{#2;#3#4}}
\newcommand\uio[2] {(u^{-1})^{#1}_{#2}}
\newcommand\uu[4]  {u_{#1}^{#2;#3#4}}
\newcommand\uuo[2] {u_{#1}^{#2}}
\newcommand\void[1]{}
\def\Vee           {^{\vee}}
\def\Xd            {{\dot X}}
\def\xx            {\_\!\_}   
\def\zet           {\ensuremath{\mathbb Z}}

\documentclass[12pt]{article}
\usepackage{latexsym,amssymb,amsfonts,bbm,epsf,color,colordvi}
\usepackage{graphics}
\usepackage[mathscr]{eucal}
\usepackage{rotating}
\begin{document}
\thispagestyle{empty}

  {~} \\[-14mm]
\begin{flushright}
{\sf math.CT/0701223}\\[1mm]{\sf KCL-MTH-07-01}\\[1mm]{\sf ZMP-HH/2007-03}\\[1mm]
{\sf January 2007} \end{flushright}
\vskip 9mm
                   ~\\[1em]~  
\begin{center}
{\large\bf THE FUSION ALGEBRA OF BIMODULE CATEGORIES}
\\[14mm]
{\large
J\"urgen Fuchs$\;^1$ \ \ \
Ingo Runkel$\;^2$ \ \ \ Christoph Schweigert$\;^{3}$} ~\\[2em]~ 
\\[6mm]
$^1\;$ Avdelning fysik, \ Karlstads Universitet\\
       Universitetsgatan 5, \ S\,--\,651\,88\, Karlstad\\[2mm]
$^2\;$ Department of Mathematics\\ King's College London, Strand\\
       GB\,--\, London WC2R 2LS\\[2mm]
$^3\;$ Organisationseinheit Mathematik, \ Universit\"at Hamburg\\
       Schwerpunkt Algebra und Zahlentheorie\\
       Bundesstra\ss e 55, \ D\,--\,20\,146\, Hamburg
\end{center}
\vskip 9mm ~\\[1em]~

\begin{quote}{\bf Abstract}\\[1mm]
We establish an algebra-isomorphism between the complexified Grothendieck 
ring \Fs\ of certain bimodule categories over a modular tensor
category and the endomorphism algebra of appropriate morphism spaces of those
bimodule categories. This provides a purely categorical proof of a conjecture
by Ostrik concerning the structure of \Fs.
\\
As a by-product we obtain a concrete expression for the structure constants of
the Grothendieck ring of the bimodule category in terms of endomorphisms of 
the tensor unit of the underlying modular tensor category.
\end{quote} 
                   \vfill \newpage  

\bigskip

\subsection*{Introduction}

The Grothendieck ring \kc\ of a semisimple monoidal category \C\ encodes a 
considerable amount of information about the structure of \C. If \C\ is 
braided, so that \kc\ is commutative, then upon complexification to 
$\kc\,{\otimes_\zet}\,\complex$ almost all of this information gets lost: what
is left is just the number of isomorphism classes of simple objects. In 
contrast, in the non-braided (but still semisimple) case the complexified 
Grothendieck ring is no longer necessarily commutative and thus contains as 
additional information the dimensions     
of its simple direct summands.

The following statement, which is a refined version of an assertion made in 
\cite{ostr} as Claim 5.3, determines these dimensions for a particularly
interesting class of categories:

\medskip

     \noindent
{\sc Theorem {\rm O}}.\\ {\it
Let \C\ be a modular tensor category, \CM\ a semisimple \ssf\ indecomposable
module category over \C, and \Cs\ the category of module endofunctors of \CM. 
Then there is an isomorphism of \complex-algebras
  \be
  \Fs \,\cong\, \bigoplus_{i,j\in I} \End_\complex\big(\HomCs(\ap(U_i),\am(U_j))
  \big)
  \labl{theiso}
between the complexified Grothendieck ring $\Fs\eq\kcs\,{\otimes_\zet}\,\complex$
and the endomorphism algebra of the specified space of morphisms in \Cs. 
}

\smallskip\noindent
Here $I$ is the (finite) set of isomorphism classes of simple objects of \C, 
$U_i$ and $U_j$ are representatives of the classes $i,j\iN I$, and $\apm$ are 
the braided-induction functors from \C\ to \Cs. 
For more details about the concepts appearing in the Theorem see section 1 
below, e.g.\ the tensor  
functors $\apm$ are given in formulas \erf{a1} and \erf{a2}.

\medskip

In \cite{ostr} the assertion of Theorem O was formulated with the help of the 
integers
  \be
  z_{i,j} := \dimc\HomCs(\ap(U_i),\am(U_j)) \,,
  \labl{zij}
in terms of which it states that \Fs\ is isomorphic to the direct sum 
$\bigoplus_{\!i,j}\mathrm{Mat}_{z_{i,j}}$ of full matrix algebras of sizes
$z_{i,j}$, $i,j\iN I$. In this form the statement had been established 
previously for the particular case that the modular tensor category \C\ is a 
category of endomorphisms of a type III factor (Theorem 6.8 of \cite{boek}), a
result which directly motivated the formulation of the statement in \cite{ostr}.
Indeed, in \cite{ostr} the additional assumption is made that the quantum 
dimension of any nonzero object of \C\ is positive, a property that is
automatically fulfilled for the categories arising in the framework of 
\cite{boek}, but is
violated in other categories (e.g.\ those relevant for the so-called 
non-unitary minimal Virasoro models) of interest in physical applications.

\medskip
 
In this note we derive the statement
in the form of Theorem O, where this positivity requirement is replaced by the
condition that \CM\ is \ssf. This property, to be explained in detail further
below, is satisfied in particular in the situation studied in \cite{boek}.
We present our proof in section 2. As an additional benefit, it provides a 
concrete expression for the structure constants of the Grothendieck ring of \Cs\
in terms of certain endomorphisms of the tensor unit of \C, see formulas 
\erf{Ns=ddd} and \erf{pic_navf_1}. Various ingredients needed in the proof are 
supplied in section 1. In section 3 we outline the particularities of the case
that \C\ comes from endomorphisms of a factor \cite{boek}, and describe 
further relations between Theorem O and structures arising in quantum field 
theory; this latter part is, necessarily, not self-contained.

      \vskip1em ~


\section{Bimodule categories and Frobenius algebras}\label{sec:bi-fro}

We start by collecting some pertinent information about the quantities used 
in the formulation of Theorem O.

\subsubsection*{Modular tensor categories}

A {\em modular tensor category\/} \C\ in the sense of Theorem O is a 
semisimple \complex-linear abelian ribbon category with simple tensor unit,
having a finite number of isomorphism classes of simple objects and obeying a 
certain nondegeneracy condition. 

Let us explain these qualifications in more detail. A {\em ribbon\/} (or 
tortile, or balanced rigid brai\-ded) category is a rigid braided monoidal 
category with a ribbon twist, i.e.\ \C\ is endowed with a tensor product 
bifunctor $\otimes$ from $\C\Times\C$ to \C, with tensor unit $\one$, and 
there are families of (right-)duality morphisms $b_U \iN\Hom(\one,U\Oti U\Vee)$,
$d_U \iN \Hom(U\Vee\Oti U,\one)$, of braiding isomorphisms $c_{U,V}
\iN\Hom(U\Oti V,V\Oti U)$, and of twist isomorphisms $\theta_U\iN\End(U)$
($U,V\iN\obj(\C)$) satisfying relations analogous to ribbons in three-space
\cite{joSt6}. (Details can e.g.\ be found in section 2.1 of \cite{fuRs4};
the category of ribbons indeed enjoys a universal property   
for ribbon categories, see e.g.\ 
chapter XIV.5.1 of \cite{KAss}.) A ribbon category is in particular 
{\em sovereign\/}, i.e.\ besides the right duality there is also a left duality,
with evaluation and coevaluation morphisms $\tilde b_U\iN\Hom(\one,{{}^\vee\!}U
\Oti U)$ and $\tilde d_U\iN \Hom(U\Oti{{}^\vee\!}U,\one)$, such that the two 
duality functors coincide, i.e.\ ${{}^\vee\!}U\eq U^\vee$ and 
${{}^\vee\!}\!f\eq f^\vee \iN \Hom(V^\vee,U^\vee)$ 
for all objects $U$ of \C\ and all morphisms $f\iN\Hom(U,V)$, $U,V\iN\obj(\C)$.

Denoting by $I$ the finite set of labels for the isomorphism classes of simple 
objects of \C\ and by $U_i$ representatives for those classes, the nondegeneracy
condition on \C\ is that the $I{\times}I$-matrix $s$ with entries
  \be\bearll
  s_{i,j} \!\!& := {\rm tr}(c_{U_i,U_j}^{}c_{U_j,U_i}^{})
  \\{}\\[-.6em]&
  \,\equiv (d_{U_j}^{}\oti\tilde d_{U_i}^{}) \cir
  [\id_{U_j^\vee}\oti(c_{U_i,U_j}^{}{\circ}\,c_{U_j,U_i}^{})\oti\id_{U_i^\vee}]
  \cir (\tilde b_{U_j}^{}\oti b_{U_i}^{}) \,\iN\End(\one) \,, \quad~ i,j\in I \,,
  \eear\labl{sij}
is invertible. We will identify $\End(\one)\eq\complex\,\id_\one$ with \complex,
so that $\id_\one\eq1\iN\complex$ and also the morphisms $s_{i,j}$ are just 
complex numbers, and we agree on $I\,{\ni}\,0$ and $U_0\eq\one$. 
In a sovereign category one has $s_{i,j}\eq s_{j,i}$. The (quantum) dimension 
of an object $U$ of a sovereign category is the trace over its 
identity endomorphism, $\dim(U)\eq\mathrm{tr}(\id_U)
\,{:=}\,d_U\cir(\id_U\oti\id_{U\Vee})\cir\tilde b_U$, 
in particular $\dim(U_i) \eq s_{0,i}$.


\subsubsection*{Module and bimodule categories}

The notion of a module category is a categorification of the one of a module
over a ring. That is, an abelian category \CM\ is a (right) {\em module 
category\/} over a monoidal category \C\ iff there exists an exact bifunctor
  \be
  \boxtimes:\quad \CM\Times\C \to \CM 
  \ee
together with functorial associativity and unit isomorphisms $(M\oxi U)\oxi V 
\Cong M\oxi(U\oti V)$ and $M\oxi\one \Cong \one$ for $M\iN\obj(\CM)$, 
$U,V\iN\obj(\C)$, which satisfy appropriate pentagon and triangle coherence 
identities. The latter involve also the corresponding coherence isomorphisms of
\C\ and are analogous to the identities that by definition of $\otimes$ are 
obeyed by the coherence isomorphisms of \C\ alone. (Any tensor category \C\ is 
a module category over itself, much like as a ring is a module 
over itself.) A module category is called {\em indecomposable\/} iff it is not 
equivalent to a direct sum of two nontrivial module categories.

A {\em module functor\/} between two module categories $\CM_1$ and $\CM_2$ over 
the same monoidal category \C\ is a functor $F{:}\ \CM_1\To\CM_2$ together with 
functorial morphisms $F(M\,{\boxtimes_1}\,U) \To F(M)\,{\boxtimes_2}\,U$ for
$M\iN\obj(\CM_1)$ and $U\iN\obj(\C)$, obeying again appropriate
pentagon and triangle identities. For details see section 2.3 of \cite{ostr}.
The composition of two module functors is again a module functor. Thus the 
category $\Cs\,{:=}\,\mathcal Fun_\C(\CM,\CM)$ of module endofunctors of a right
module category \CM\ over a monoidal category \C\ is monoidal, and the action of
these endofunctors on \CM\ turns \CM\ into a left module category over \Cs.

If \C\ has a braiding $c$ then two functors $\apm$ from \C\ to \Cs\ are of
particular interest. They are defined by assigning to $U\iN\obj(\C)$ the
endofunctors
  \be
  \ap(U) = \am(U) := \xx \oxi U
  \labl{a1}
of \CM. In order that $\apm$ are indeed functors to \Cs, we have
to make $\apm(U)$ into module functors, i.e.\ to specify morphisms
$\apm(U)(M\Oxi V) \To \big(\apm(U)(M)\big) \oxi V$, i.e.
  \be
  M \oxi (V\oti U) \cong (M\oxi V) \oxi U 
  ~\longrightarrow~ (M\oxi U) \oxi V \cong M\oxi (U\oti V)
  \labl{a2}
for $M\iN\obj(\CM)$ and $U,V\iN\obj(\C)$; this is achieved by using the 
morphisms\,%
  \footnote{~suppressing, for brevity, the mixed associativity morphisms}
$\id_M \oti c_{V,U}$ for $\ap$ and $\id_M \oti c_{U,V}^{-1}$ for $\am$,
respectively. 

The functors $\apm$ defined this way are monoidal; we call them 
{\em braided-induction\/} functors. These functors, originally referred to as 
{\em alpha induction\/} functors, were implicitly introduced \cite{lore} and 
heavily used \cite{xu3,boev,boek,boev5} in the context of subfactors, and 
were interpreted as monoidal functors with values in a category of module 
endofunctors in \cite{ostr}. 

Further, by setting
  \be
  (U,X,V) \,\mapsto\, \beta^{\Eps\Eps'}(U,V)
  := \alpha^{\Eps}(U) \cir X \cir \alpha^{\Eps'}(V)
  \ee
for $U,V\iN\obj(\C)$ and $X\iN\obj(\Cs)$, for any choice of signs $\Eps,\Eps'
\iN\{\pm\}$ one obtains a functor $\beta^{\Eps\Eps'}{:}\ \C\Times\Cs\Times\C
\To\Cs$; it can be complemented by two sets of associativity constraints
(for the left and right action of \C), both obeying separately a pentagon
constraint and an additional mixed constraint that expresses the
commutativity of the left and right action of \C. This
turns \Cs\ into a {\em bimodule category\/} over \C.


\subsubsection*{Frobenius algebras}

Given a monoidal category \C, the structure of a module
category over \C\ on an abelian category \CM\ is equivalent to
a monoidal functor from \C\ to the category of endofunctors of \CM.
If \C\ is semisimple rigid monoidal, with simple tensor unit and with
a finite number of isomorphism classes of simple objects, and \CM\ is 
indecomposable and semisimple, then it follows \cite[Theorem\,3.1]{ostr}
that \CM\ is equivalent to the category \CA\ of (left) $A$-modules in \C\ for
some algebra $A$ in \C. By a similar reasoning \Cs\ is then equivalent to the 
category \CAA\ of $A$-bimodules in \C. The monoidal product of \CAA\ is the 
tensor product over $A$. The algebra $A$ such that $\CA\,{\simeq}\,\CM$ is 
determined uniquely up to Morita equivalence; it can be constructed as
the internal End $\underline{\End}(M)$ for any $M\iN\obj(\CM)$ \cite{ostr}.

Recall that a (unital, associative) algebra $A\eq (A,m,\eta)$ in a (strict) 
monoidal category \C\ consists of an object $A\iN\obj(\C)$ and morphisms 
$m \iN \Hom(A\oti A,A)$ and $\eta \iN \Hom(\one,A)$ satisfying 
$m\cir(m\oti\ida) \eq m\cir(\ida\oti m)$ and
$m\cir(\eta\oti\ida) \eq \ida \eq m\cir(\ida\oti\eta)$.
If \CM\ is indecomposable then the category \Cs\ has simple tensor unit (it is 
thus a fusion category in the sense of \cite{etno}). In terms of the algebra 
$A$, this property means that $A$ is simple as an object of the category 
\CAA\ of $A$-bimodules in \C; such algebras $A$ are called {\em simple\/}.

A left $A$-module is a pair $M\eq(\dot M,\rho)$ consisting of an object 
$\dot M\iN\obj (\C)$ and of a morphism $\rho\iN\Hom(A\oti\dot M,\dot M)$ that 
satisfies $\rho\cir(\ida\oti\rho) \eq 
\rho\cir(m\oti\idM)$ and $\rho\cir(\eta\oti\idM) \eq \idM$. A right $A$-module
$(\dot M,\varrho)$, $\varrho\iN\Hom(\dot M,\dot M\oti A)$, is defined 
analogously, and an $A$-bimodule is a triple $X\eq(\Xd,\rhol,\rhor)$ such that 
$(\Xd,\rhol)$ is a left $A$-module, $(\Xd,\rhor)$ is a right $A$-module
and the left and right $A$-actions commute. The morphism space $\HomA(M,N)$ 
in \CA\ consists of those morphisms in $\Hom(\dot M,\dot N)$ which commute with 
the left $A$-action, and an analogous property characterizes the morphism space 
$\HomAA(X,Y)$ in \CAA.

\medskip

Of interest to us is a particular class of algebras in \C, the symmetric special
Frobenius algebras. A {\em Frobenius algebra\/} in a monoidal category
\C\ is a quintuple $A\eq(A,m,\eta,\Delta,\eps)$ such that $(A,m,\eta)$ is
an algebra, $(A,\Delta,\eps)$ is a coalgebra, and the coproduct $\Delta$ is
a morphism of $A$-bimodules.\,%
  \footnote{~In the classical case of Frobenius algebras in the category of
  vector spaces over a field, the Frobenius property can be formulated in
  several other equivalent ways. The one given here does not require any
  further structure on \C\ beyond monoidality. Also note that neither $\Delta$
  nor the counit $\eps$ is required to be an algebra morphism.}
A Frobenius algebra $A$ in a sovereign monoidal category is {\em symmetric\/} 
iff the morphism $(d_A\oti\ida) \cir [\idav \oti (\Delta\cir\eta\cir\eps\cir m)]
\cir (\tilde b_A\oti\ida) \iN\End(A)$
(which for any Frobenius algebra $A$ is an algebra automorphism) equals $\ida$. 
A Frobenius algebra $A$ is {\em special\/} iff $\Delta$ is a
right-inverse of the product $m$ and the counit $\eps$ is a left-inverse of 
the unit $\eta$, up to nonzero scalars. For a special Frobenius algebra one
has $\dim(A)\nE 0$, and one can normalize the counit 
in such a way that $m\cir\Delta\eq\ida$ and $\eta\cir\eps\eq\dim(A)\,\id_\one$;
below we assume that this normalization has been chosen.

The structure and \rep\ theory of symmetric special Frobenius algebras have been
studied e.g.\ in \cite{kios,fuRs4,ffrs,fuRs11}.
The braided-induction functors $\apm$ which exist when \C\ has a braiding can be
described in terms of the algebra $A$ as the functors $\apm_{\!A}{:}\ \C\To  
\CAA\,{\simeq}\,\Cs$ that associate to $U\iN\obj(\C)$ the following 
$A$-bi\-modules $\apm_{\!A}(U)$:      
the underlying object is $A\oti U$, the left module 
structure is the one of an induced left module, and the right module structure 
is given by the one of an induced right module composed with a braiding 
$c_{U,A}^{~-1}$ for $\ap_{\!A}$ and $c_{\!A,U}^{}$ for $\am_{\!A}$,  
respectively \cite{ostr,ffrs}. 
In terms of $A$ the numbers $z_{ij}$ in \erf{zij} are given by
$z_{ij} \eq z(A)_{ij}$ with       
  \be
  z(A)_{i,j} := \dimc\HomAA(\ap_{\!A}(U_i),\am_{\!A}(U_j)) \,.
  \labl{zaij}
Let us also mention that for modular $\C$ every $A$-bimodule 
can be obtained as a retract of a tensor product
$\ap_{\!A}(U)\,{\otimes_{\!A}^{}}\,\am_{\!A}(V)$ of a suitable pair
of $\ap_{\!A}$- and $\am_{\!A}$-induced      
bimodules (see the conjecture \cite[Claim\,5.2]{ostr} and its proof in 
\cite{ffrs5}), and  that the matrix $z(A)$      
is a permutation matrix iff the functors $\apm_{\!A}$     
are monoidal equivalences, i.e.\ iff $A$ is Azumaya \cite{fuRs11}.


\subsubsection*{\Ssfy}

To cover the terms used in the formulation of Theorem O we need to introduce
one further notion, the one of \ssfy. This is done in the following

     \noindent
{\sc Definition}.\\
(i)~\,{\it An algebra $A$ in a sovereign tensor category is called {\em \ssf\/}
iff the morphism
  \be
  \big[ (\eps_\natural \cir m) \oti \idav \big] \circ (\ida \oti b_A)
  ~\in \Hom(A,A^\vee_{}) 
  \labl{Phi1}
with $\eps_\natural \,{:=}\, d_A \cir (\idav \oti m) \cir (\tilde b_A \oti \ida)
\iN \Hom(A,\one)$ is an isomorphism.}
\\[2pt]
(ii)~{\it A semisimple indecomposable module category \CM\ over a semisimple 
sovereign tensor category \C\ that has simple tensor unit and a finite number 
of isomorphism classes of simple objects is called {\em \ssf\/} iff
there exists a \ssf\ algebra $A$ in \C\ such that $\CM\,{\simeq}\,\CA$.}

\medskip

It is not difficult to see that the property of an algebra to be \ssf\ is 
preserved under Morita equivalence. But it is at present not evident to us how 
to give a definition of \ssfy\ for module categories which does not make direct
reference to the corresponding (Morita class of) algebra(s); we plan to come 
back to this problem elsewhere. That such a formulation must exist is actually 
the motivation for introducing this terminology. In contrast, for algebras 
\ssfy\ is not a new concept, owing to

\smallskip

     \noindent
{\sc Lemma 1} \\ {\it An algebra in a sovereign tensor category is \ssf\
iff it is symmetric special Frobenius.}

\smallskip

     \noindent Proof.
\\
That a \ssf\ algebra in a sovereign tensor category is symmetric special 
Frobenius has been shown in lemma 3.12 of \cite{fuRs4}. The converse holds
because for a symmetric special Frobenius algebra one has $\eps_\natural 
\eq \eps$, and as a consequence
$(d_A\oti\ida)\cir(\idav\oti(\Delta\cir\eta)) \iN \Hom(A^\vee_{},A)$
is inverse to the morphism \erf{Phi1}.
\qed


\subsubsection*{Fusion rules}

We denote by $[U]$ the isomorphism class of an object $U$.
If an abelian category \C\ is rigid monoidal, then the tensor product bifunctor 
is exact, so that the Grothendieck group \kc\ has a natural ring structure given
by $[U]\Star[V]\,{:=}\,[U\oti V]$. A distinguished basis of \kc\ is given by the
isomorphism classes of simple objects of \C; in this basis the structure 
constants are non-negative integers. If \CM\ is a module category over \C, then 
the Grothendieck group $K_0(\CM)$ is naturally a \kc-module, in fact a based 
module over the based ring \kc.

{}From now on, unless noted otherwise, \C\ will stand for a modular tensor 
category, \CM\ for a semisimple \ssf\ indecomposable module category over \C, and
$\Cs\eq\mathcal Fun_\C(\CM,\CM)$ will be regarded as a bimodule category over 
\C\ via the functor $\bpm$. Then the category \Cs\ is semisimple 
\complex-linear abelian rigid monoidal and has \findim\ morphism spaces and 
a finite number of isomorphism classes of simple objects. 
For simplicity of notation we also tacitly take \C\ 
to be strict monoidal; in the non-strict case the relevant coherence 
isomorphisms have to be inserted at appropriate places, 
but all statements about Grothendieck rings remain unaltered.

Analogously as \kc\ also the Grothendieck group \kcs\ has a natural ring 
structure. However, \Cs\ is not, in general, braided, and hence, unlike \kc, 
the ring \kcs\ is in general not commutative.
The complexified Grothendieck ring 
  \be
  \Fs := \kcs \otimes_\zet \complex
  \ee
of \Cs\ is a \findim\ semisimple associative \complex-algebra, and hence a 
direct sum of full matrix algebras, $\Fs\Cong\bigoplus_{p\in P}\Mat_{n_p}$ for 
some finite index set $P$, with $\Mat_n$ the algebra of $n\Times n$\,-matrices 
with complex entries. We denote the product in \Fs\ by $\Star$. The algebra \Fs\
has a standard basis $\{ \e p\alpha\beta \}$, $\Fs \Cong \bigoplus_{p\in P}
\bigoplus_{\alpha,\beta=1}^{n_p} \complex\e p\alpha\beta$, with products
  $
  \e p\alpha\beta \Star \e q\gamma\delta
  \eq \delta_{p,q}\,\delta_{\beta,\gamma}\, \e p\alpha\delta
  $.
The elements $e_p\,{:=}\,\sum_{\alpha=1}^{n_p}\!\e p\alpha\alpha$ are the 
primitive idempotents projecting onto the simple summands $\Mat_{n_p}$ of \Fs.
Given any 
other basis $\{ x_\ka \,|\, \ka\iN\Is \}$, there is a basis transformation
  \be 
  x_\ka = \sum_{p\in P}\sum_{\alpha,\beta=1}^{n_p}
  \uu \ka p\alpha\beta \e p\alpha\beta \,, \qquad
  \e pab = \sum_{\ka\in K} \ui \ka p\alpha\beta \, x_\ka \,,
  \labl{uu}
satisfying
  $
  \sum_{p,\alpha,\beta} \ui \kap p\alpha\beta\, \uu \ka p\alpha\beta
  \eq \delta_{\ka,\kap} $ and $
  \sum_\ka \uu \ka q\alpha\beta \, \ui \ka p\gamma\delta
  \eq \delta_{p,q}\,\delta_{\alpha,\gamma}\,\delta_{\beta,\delta} \,,
  $
so that the structure constants of \Fs\ in the basis $\{ x_\ka \}$ can
be written as
  \be
  \Ns \ka\kap\kapp = \sum_{p\in P}\sum_{\alpha,\beta,\gamma=1}^{n_p}
  \uu \ka p\alpha\beta\, \uu \kap p\beta\gamma\, \ui \kapp p\alpha\gamma \,.
  \labl{Ns=uuu}
For the particular case that $\{ x_\ka \}$ is the distinguished basis 
of the underlying ring given by the isomorphism classes $[X_\ka]$ of simple 
objects of \Cs, it has become customary in various contexts to refer to the 
structure constants $\Ns \ka\kap\kapp$
of \Fs\ (or also to the algebra \Fs, or to the Grothendieck ring itself) as 
the {\em fusion rules\/} of the category under study. The formula \erf{Ns=uuu} 
then constitutes a block-diagonalization of the fusion rules of \Cs.
Also note that if \Fs\ is commutative, then \erf{Ns=uuu} reduces to
$\Ns\ka\kap\kapp \eq \sum_{p\in P} \uuo\ka p\,\uuo \kap p\,\uio\kapp p$, which
is sometimes referred to as the Verlinde formula.

\medskip

In terms of these ingredients, the assertion of Theorem O amounts to the claim
that the index set $P$ and the dimensions $n_p^2$ of the simple summands of
\Fs\ are given by
  \be
  P = \{ (i,j)\iN I{\times}I \,|\, z_{i,j} \nE 0 \} \qquad{\rm and}\qquad
  n_{(i,j)} = z_{i,j} 
  \labl{P-claim}
with the integers $z_{i,j}$ the dimensions defined in formula \erf{zij}.

We will establish the equalities \erf{P-claim} by finding an explicit 
expression for the matrix elements 
$\uu \ka p\alpha\beta \,{\equiv}\, \uu \ka {(i,j)}\alpha\beta$ of the basis
transformation \erf{uu} for the case that the basis $\{x_\ka \,|\, \ka\iN\Is \}$
is the distinguished basis   
of isomorphism classes $[X_\ka]$ of the simple objects.


\section{The structure of the algebra \boldmath{\Fs}}

The considerations above show in particular that a nondegenerate
semisimple indecomposable module category $\CM$ over a modular tensor
category \C\ is equivalent, as a module category, to the category
\CA\ for an appropriate simple symmetric special Frobenius algebra $A$ in \C. 
Moreover, the category \Cs\ of module endofunctors is equivalent to \CAA\ 
as a monoidal category (and as a bimodule category over \C\ and left 
module category over \CM). We may thus restate Theorem O as

     \medskip\noindent
{\sc Theorem ${\rm O'}$}.\\ {\it
Let $\C$ be a modular tensor category and $A$ a simple symmetric special 
Frobenius algebra in $\C$. Then the complexified Grothendieck ring
$\Fs\eq\kcaa\,{\otimes_\zet}\,\complex$ and the endomorphism algebra
  \be
  \Es
  := \bigoplus_{i,j\in I}\End_\complex\big(\HomAA(\ap(U_i),\am(U_j))\,.
  \labl{Es}
are isomorphic as \complex-algebras.}

\medskip

We will prove Theorem ${\rm O}'$ by constructing an 
isomorphism from $\Fs$ to $\Es$.

\subsubsection*{A map \boldmath{$\Phi$} from \boldmath{\Fs} to \boldmath{\Es}}

To simplify some of the expressions appearing below, we replace $j\iN I$ in 
\erf{zij} by $\overline\jmath$, defined as the unique label in $I$ such that
$U_{\overline\jmath}^{}\Cong U_j^\vee$,
and instead of the spaces $\HomAA(\ap(U),\am(V^\vee))$
we prefer to work with the isomorphic spaces $\HomAA(U\otip A\otim V,A)$, where 
for any $A$-bimodule $X\eq(\Xd,\rhol,\rhor)$ the $A$-bi\-mo\-dules $U\otip X$
and $X\otim V$ are defined as $U\otip X \,{:=}\, (U \Oti \Xd , (\id_U\Oti\rhol)
      \,{\circ}$\linebreak[0]$
(c_{U,A}^{~-1}\Oti\id_{\Xd}) , \id_U\Oti\rhor)$ and as
$X\otim V \,{:=}\, (\Xd \Oti V , \rhol\Oti\id_V, (\rhor\Oti\id_V) \cir 
(\id_{\Xd}\Oti c_{\!A,V}^{\,-1}))$, respectively.
Accordingly instead of with \erf{Es} we work with the isomorphic algebra
  \be
  \bigoplus_{i,j\in I}\End_\complex\big(\HomAA(U_i\otip A\otim U_j,A)\big)
  \cong \Es
  \ee
which by abuse of notation we still denote by \Es.

For any $U,V\iN\obj(\C)$ and $X\iN\obj(\CAA)$ we denote by $\D XUV$ the mapping
  \begin{eqnarray}
  && \hspace*{-1.3em} \bearl
  \varphi \in \HomAA(U\otip A\otim V,A)  
  \\{}\\[-.8em] \hspace*{.1em}
  \longmapsto~ \D XUV(\varphi) :=
  \big( \ida \oti \tilde d_\Xd \big)
  \circ \Big[ \Big( \big[ \ida \oti \rhol \oti (\eps\cir\varphi) \big]    
                    \circ \big[ (\Delta\cir\eta) \oti c_{U,\Xd}
                    \oti \id_A \oti \id_V \big] \Big) \oti \id_{\Xd^\vee} \Big]
  \\{}\\[-.8em] \hspace*{7.8em}
  \circ\, \Big[ \id_U \oti
                \Big( \big[ (\rhor \cir (\rhol \oti \ida)) \oti \ida \big]
                      \cir \big[ \ida \oti \id_\Xd \oti (\Delta\cir\eta) \big] \Big)
                \oti \id_V \oti \id_{\Xd^\vee} \Big]
  \\{}\\[-.8em] \hspace*{7.8em}
  \circ\, \Big[ \id_U \oti \ida \oti
                \Big( \big[ c_{\Xd,V}^{-1} \oti \id_{\Xd^\vee} \big]
                      \cir \big[ \id_V \oti b_\Xd \big] \Big) \Big]
  \eear\nonumber\\[-1.6em]~
  \label{DXUV}\end{eqnarray}
One checks that $\D XUV(\varphi)$, which by construction is a morphism in
$\Hom(U\oti A\oti V,A)$, is actually again in the subspace
$\HomAA(U\otip A\otim V,A)$, so that $\D XUV$ is an endomorphism of the vector 
space $\HomAA(U\otip A\otim V,A)$. 
Further, we consider the linear map 
  \be
  \Phi: \quad \Fs\eq\kcaa\,{\otimes_\zet}\,\complex \,\to\, \Es
  \ee 
defined by
  \be
  \Phi([X]) := \D X{}{} = \bigoplus_{i,j\in I} \D X{U_i}{U_j} .
  \labl{Phi}

It must be admitted that the expression \erf{DXUV} for the map $\D XUV$ is not 
exceedingly transparent. However, with the help of the graphical notation for 
monoidal categories as described e.g.\ in \cite{joSt6,MAji,KAss},
it can easily be visualized. Indeed, using
  \eqpic{pic_navf_a} {400} {31} {
  \put(20,0)     {\Includeourbeautifulpicture{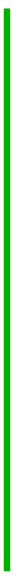}} 
  \put(-23,30)   {$ \id_U^{}\,= $}
  \put(18.0,-8.8){\sse$ U $}
  \put(18.5,65.5){\sse$ U $}
\put(87,0){
  \put(20,0)     {\Includeourbeautifulpicture{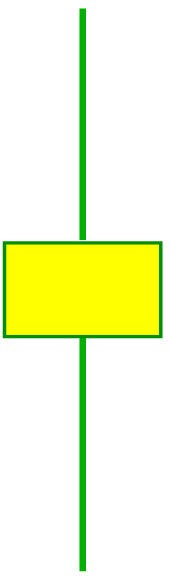}}
  \put(-16,30)   {$ f~= $}
  \put(25.3,-8.8){\sse$ U $}
  \put(26.2,65.5){\sse$ V $}
  \put(26.6,30.2){\sse$ f $}
}
\put(204,0){
  \put(20,0)     {\Includeourbeautifulpicture{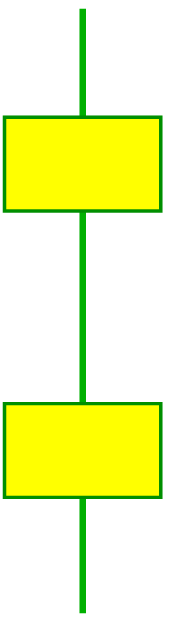}}
  \put(-31,30)   {$ g\cir f~= $}
  \put(25.1,69.5){\sse$ W $}
  \put(25.3,-8.8){\sse$ U $}
  \put(26.6,17.6){\sse$ f $}
  \put(26.6,49.5){\sse$ g $}
  \put(30.5,32.9){\sse$ V $}
}
\put(330,0){
  \put(20,0)     {\Includeourbeautifulpicture{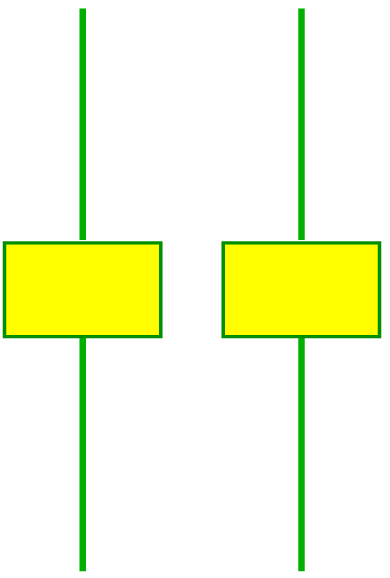}}
  \put(-40,30)   {$ f\oti f'~= $}
  \put(25.3,-8.8){\sse$ U $}
  \put(26.2,65.5){\sse$ V $}
  \put(26.6,29.4){\sse$ f $}
  \put(49.1,-8.8){\sse$ U' $}
  \put(50.0,65.5){\sse$ V' $}
  \put(49.6,29.4){\sse$ f' $}
}
  }
for identity morphisms, general morphisms $f \iN \Hom(U,V)$, 
and for composition and tensor product of morphisms of \C,
  \eqpic{pic_navf_b} {380} {46} {
\put(0,22){
  \put(10,0)     {\Includeourbeautifulpicture{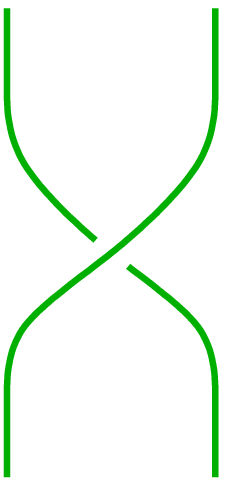}} 
  \put(-39,23)   {$ c_{U,V}^{}~= $}
  \put(7.2,-8.8) {\sse$ U $}
  \put(7.8,56.3) {\sse$ V $}
  \put(30.2,-8.8){\sse$ V $}
  \put(31.9,56.3){\sse$ U $}
}
\put(122,22){
  \put(10,0)     {\Includeourbeautifulpicture{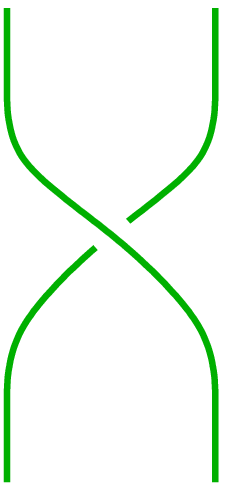}} 
  \put(-40,23)   {$ c_{U,V}^{-1}~= $}
  \put(7.2,-8.8) {\sse$ V $}
  \put(7.8,56.3) {\sse$ U $}
  \put(30.2,-8.8){\sse$ U $}
  \put(31.9,56.3){\sse$ V $}
}
\put(239,60){
  \put(10,0)     {\Includeourbeautifulpicture{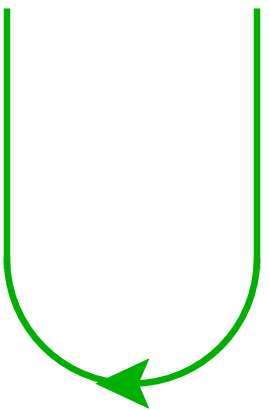}} 
  \put(-32,22)   {$ b_U~= $}
  \put(8.3,47.3) {\sse$ U $}
  \put(33.9,47.3){\sse$ U^\vee $}
}
\put(239,0){
  \put(10,0)     {\Includeourbeautifulpicture{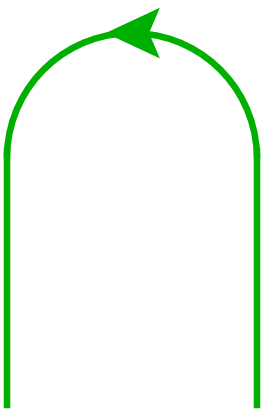}} 
  \put(-32,16)   {$ d_U~= $}
  \put(6.2,-8.2) {\sse$ U^\vee $}
  \put(34.8,-8.2){\sse$ U $}
}
\put(360,60){
  \put(10,0)     {\Includeourbeautifulpicture{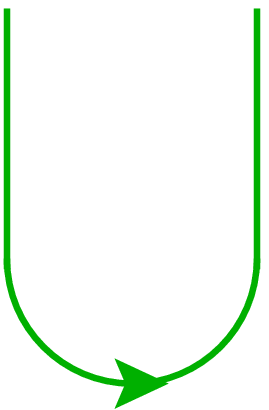}} 
  \put(-32,22)   {$ \tilde b_{U}^{}~= $}
  \put(-2.4,47.3){\sse$ {\phantom U}^\vee_{}\!{U} $}
  \put(34.5,47.3){\sse$ U $}
}
\put(360,0){
  \put(10,0)     {\Includeourbeautifulpicture{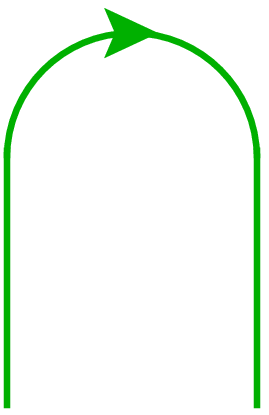}} 
  \put(-32,16)   {$ \tilde d_{U}^{}~= $}
  \put(7.8,-8.2) {\sse$ U $}
  \put(24.9,-8.2){\sse$ {\phantom U}^\vee_{}\!{U} $}
}
  }
for braiding and duality morphisms of \C, as well as
                   \newpage~\\[-1.5em] 
  \eqpic{pic_navf_c} {380} {47} {
\put(0,58){
  \put(10,0)     {\Includeourbeautifulpicture{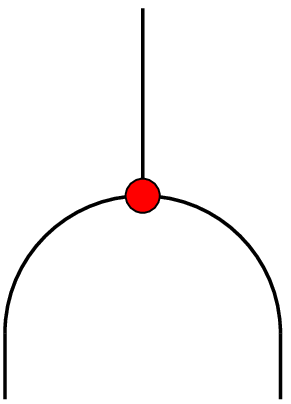}} 
  \put(-29,20)   {$ m~= $}
  \put(7.5,-8.8) {\sse$ A $}
  \put(22.5,46.5){\sse$ A $}
  \put(37.5,-8.8){\sse$ A $}
}
\put(14,-6){
  \put(10,0)     {\Includeourbeautifulpicture{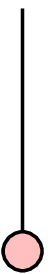}} 
  \put(-25,14)   {$ \eta~= $}
  \put(9.5,32.2) {\sse$ A $}
}
\put(120,58){
  \put(10,0)     {\Includeourbeautifulpicture{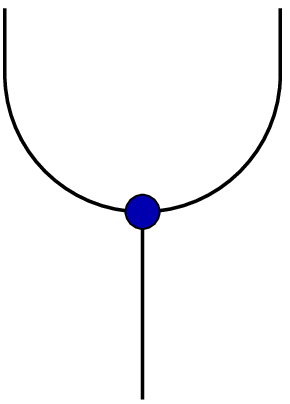}}
  \put(-29,20)   {$ \Delta~= $}
  \put(7.5,46.5) {\sse$ A $}
  \put(22.5,-8.8){\sse$ A $}
  \put(37.5,46.5){\sse$ A $}
}
\put(134,0){
  \put(10,0)     {\Includeourbeautifulpicture{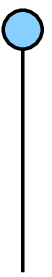}} 
  \put(-25,14)   {$ \eps~= $}
  \put(8.8,-9.6) {\sse$ A $}
}
\put(240,18){
  \put(10,0)     {\Includeourbeautifulpicture{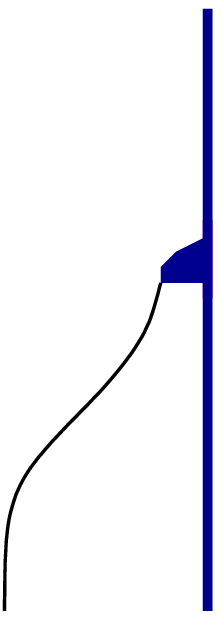}}
  \put(-26,23)   {$ \rhol~= $}
  \put(7.5,-8.8) {\sse$ A $}
  \put(29.5,-9.4){\sse$ \Xd $}
  \put(30.5,69.4){\sse$ \Xd $}
}
\put(360,18){
  \put(10,0)     {\Includeourbeautifulpicture{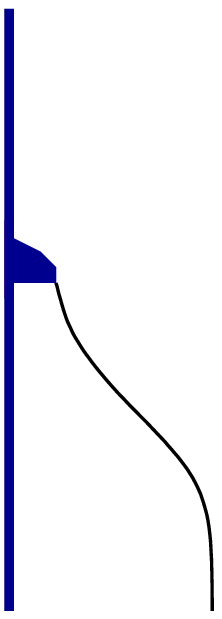}}
  \put(-31,23)   {$ \rhor~= $}
  \put(29.5,-8.8){\sse$ A $}
  \put(7.5,-9.4) {\sse$ \Xd $}
  \put(8.5,69.4) {\sse$ \Xd $}
}
  }
for the structural morphisms of the algebra $A$ and of the bimodule $X$,
the definition \erf{DXUV} amounts to 
  \eqpic{pic_navf_6} {183} {97} {
      \put(0,3){
  \put(75,0)      {\Includeourbeautifulpicture 6}
  \put(0,94)      {$ \D X UV(\varphi)~= $}
  \put(79.5,-8.8) {\sse$ U $}
  \put(89.5,-8.8) {\sse$ A $}
  \put(89.5,209)  {\sse$ A $}
  \put(100.5,-8.8){\sse$ V $}
  \put(103.8,87.5){\sse$ \Xd $}
  \put(138.5,157.8){\sse$ \varphi $}
  } }
With the help of this graphical description it is e.g.\ easy to verify
that $\D X UV$ only depends on the isomorphism class $[X]$ of the bimodule
$X$: Given any bimodule isomorphism $f\iN\HomAA(X,X')$ one may insert
$\id_X \eq f^{-1}\cir f$ anywhere in the $X$-loop and then 
`drag $f$ around the loop' and use $f\cir f^{-1}\eq\id_{X'}$, thereby
replacing the $X$-loop by an $X'$-loop.


\subsubsection*{\boldmath{$\Phi$} is an algebra morphism}

{\sc Lemma 2} \\ {\it
The map $\Phi$ defined in \erf{Phi} is a morphism of unital associative 
\complex-algebras.}

\medskip

     \noindent Proof.
\\
We need to show that $\Phi(1_{\Fs})\eq 1_{\Es}$ and
that $\Phi([X])\cir\Phi([Y]) \eq \Phi([X]\Star[Y])$ for all $X,Y\iN\obj
  $\linebreak[0]$
(\CAA)$, 
or in other words, that
  \be
  \D {\!A}{}{} = \id     \qquad{\rm and}\qquad
  \D X{}{} \circ \D Y{}{} = \D {X\Ota Y}{}{} \,,
  \labl{11}
where $\ota$ is the tensor product over $A$. The latter equality is seen as 
follows. Using the defining property of the counit, the unitality of the left
$A$-action $\rhol$, as well as the functoriality of the braiding and the fact
that the left and right $A$-actions on the bimodule $Y$ commute, one obtains
  \eqpic{pic_navf_7} {280} {134} {
  \put(110,0)    {\Includeourbeautifulpicture 7}
  \put(0,136)    {$ \D X{}{} \circ \D Y{}{}(\varphi) ~= $}
  \put(114.5,-8.8){\sse$ U $}
  \put(124.5,-8.8){\sse$ A $}
  \put(124.5,284) {\sse$ A $}
  \put(135.5,-8.8){\sse$ V $}
  \put(139.3,147) {\sse$ \Xd $}
  \put(166.5,103) {\sse$ \dot Y $}
  \put(205.5,172.3){\sse$ \varphi $}
  }
The assertion then follows immediately by the fact \cite{kios,ffrs5}
that the morphism
  \eqpic{pic_navf_5} {250} {20} {
  \put(200,0)    {\Includeourbeautifulpicture 5}
  \put(0,24)     {$ (\rhor^X \oti \rhol^Y) \circ
                    (\id_\Xd \oti (\Delta\cir\eta) \oti \id_{\dot Y}) ~= $}
  \put(196.8,-8.8){\sse$ \Xd $}
  \put(210.2,31.7){\sse$ A $}
  \put(224.5,-8.8){\sse$ \dot Y $}
  }
is the idempotent corresponding to the epimorphism that restricts $X\oti Y$ 
to $X\ota Y$.
\\[2pt]
To show also the first of the equalities \erf{11}, first note that for $X\eq A$ 
the left and right $A$-actions are just given by the product $m$ of $A$. One can
then use the various defining properties of $A$ and the fact that $\varphi\iN
\HomAA(U\otip A\otim V,A)$ intertwines the left action of $A$ to arrive at
                   \newpage~\\[-1.5em] 
  \eqpic{pic_navf_8} {220} {110} {
     \put(0,3){
  \put(100,0)    {\Includeourbeautifulpicture 8}
  \put(30,112)   {$ \D {\!A}UV(\varphi) ~= $}
  \put(96.5,-8.8) {\sse$ U $}
  \put(106.5,-8.8){\sse$ A $}
  \put(106.5,224) {\sse$ A $}
  \put(117.5,-8.8){\sse$ V $}
  \put(136.3,76.8){\sse$ A $}
  \put(148.3,169.8){\sse$ A $}
  \put(156.5,152.3){\sse$ \varphi $}
  } }
By again using the associativity, symmetry and Frobenius property of $A$, 
together with the fact that $\varphi$ is a bimodule intertwiner, one can reduce
the right hand side to $\varphi$, except for an additional $A$-loop being 
attached to the $A$-line. This $A$-loop, in turn, is removed by
invoking specialness of $A$. Thus $\D{\!A}UV$ acts as the
identity on $\HomAA(U\otip A\otim V,A)$, which establishes the claim.
\qed


\subsubsection*{\boldmath{$\Phi$} is an isomorphism}

Maps $\D X{}{}$ of the type appearing in \erf{Phi} have already been considered in
\cite{ffrs5}, where their properties were studied using the relation
\cite{TUra} between modular tensor categories and three-dimensional topological
field theory. Indeed, these maps are special cases of certain linear maps
$\Hom_{B|B}(U\otip B\otim V,B)\To\HomAA(U\otip\, Y\otim V,A)$ defined for
a pair $A$, $B$ of symmetric special Frobenius algebras and an $A$-bimodule $Y$,      
see equations (2.29) and and (4.14) of \cite{ffrs5}.

\medskip

     \noindent
{\sc Lemma 3} \\ {\it
The algebra morphism $\Phi$ defined in \erf{Phi} is an isomorphism.
}

\smallskip

     \noindent Proof.
\\
It has been shown in Proposition 2.8 of \cite{ffrs5} that if 
the equality $\D X{}{}(\varphi) \eq \D Y{}{}(\varphi)$ holds for all 
$\varphi\iN\HomAA(U_i\otip A\otim U_j,A)$ and $i,j\iN I$, then
$X$ and $Y$ are isomorphic bimodules. Hence $\Phi$ is injective.
\\
That $\Phi$ is surjective thus follows from the fact that
  \be
  \dimc(\Fs)= {\rm tr} \big( z(A)^{\rm t}_{} z(A) \big)
  \equiv \sum_{i,j\in I} \big( z(A)_{i,j} \big)^2_{} = \dimc(\Es) \,.
  \ee
Here the first equality has been established in Remark 5.19(ii) of \cite{fuRs4},
while the last equality is satisfied because, by definition of the numbers 
$z(A)_{i,j}$, 
$\dimc\big(\End_\complex(\HomAA(U_i\otip A\otim U_{\overline\jmath},A)) \big) 
       $\linebreak[0]${=}\,
{(z(A)_{i,j})}^2_{}$.
\qed

\medskip

This completes the proof of Theorem ${\rm O}'$, and thus also of Theorem 
${\rm O}$.
\qed


\subsubsection*{Block-diagonalization}

Using the fact that for simple $A$ the mapping $f\,{\mapsto}\,\eps\cir f\cir
\eta$ furnishes a canonical isomorphism $\EndAA(A)\To\End(\one)\eq\complex$, 
one sees that with respect to the bases of $\Fs$ and $\Es$
considered in section \ref{sec:bi-fro} the 
linear map $\Phi$ \erf{Phi} is given by the matrix \rmd\ with entries
  \be
  \dd \ka ij\alpha\beta
  := \eps \circ \D{X_\ka}{U_i}{U_j}(h^{(ij)}_\beta) \circ \ol h^{(ij)}_\alpha
  \cir \eta ~\in \End(\one)\eq\complex \,,
  \labl{dd}
where $\ka\iN\Is$, $i,j\iN I$ and $h^{(ij)}_\beta\iN\Ga ij$, $\ol h^{(ij)}
_\alpha\iN\Gab ij$, with $\Ga ij$ a basis of $\HomAA(U_i\otip A\otim U_j,A)$ and
$\Gab ij$ a dual basis of $\HomAA(A,U_i\otip A\otim U_j)\Cong \HomAA{(U_i\otip A
\otim U_j,A)}^*$. Since $\Phi$ is an isomorphism, \rmd\ is invertible.

Now $\Phi(x_\kappa)$ acts on the elements of a basis $\bigcup_{i,j\in I}\Ga ij$
of $\bigoplus_{i,j \in I}\HomAA(U_i\otip A\otim U_j,A)$         
as $h^{(ij)}_\alpha\,{\mapsto}\,\sum_\beta\dd \ka ij\beta \alpha
h^{(ij)}_\beta$, while the elements $\Phi(\ei kl\gamma\delta)$ are by definition
the matrix units of \Es\ with respect to the basis $\{h^{(ij)}_\alpha\}$,
mapping $h^{(ij)}_\alpha$ to $\delta_{ik}\delta_{jl}\delta_{\alpha\delta}\,
h^{(ij)}_\gamma$, and hence we have $\Phi(x_\kappa) \eq \sum_{i,j,\alpha,\beta}
\dd\ka ij\alpha\beta
        $\linebreak[0]$     
\Phi(\ei ij\alpha\beta)$. Thus \rmd\ is precisely the 
basis transformation matrix $u$ appearing in formula \erf{uu}. As a consequence,
according to \erf{Ns=uuu} the fusion rules of \CAA\ can be written as
  \be
  \Ns \ka\kap\kapp
  = \sum_{i,j\in I}\sum_{\alpha,\beta,\gamma=1}^{z_{i,\overline\jmath}^{}}
  \dd \ka ij\alpha\beta\, \dd \kap ij\beta\gamma\, \ddi \kapp ij\alpha\gamma \!.
  \labl{Ns=ddd}
Pictorially, the entries of the block-diagonalization matrix are given by
  \eqpic{pic_navf_1} {380} {103} {
  \put(75,0)       {\Includeourbeautifulpicture 1}
  \put(245,0)      {\Includeourbeautifulpicture 4}
  \put(0,104)      {$ \dd \ka ij\alpha\beta~= $}
  \put(65.7,77)    {\sse$ U_i $}
  \put(101.5,173)  {\sse$ \Xd_\ka^{} $}
  \put(131.2,18.2) {\sse$ \ol\alpha $}
  \put(131.2,181.2){\sse$ \beta $}
  \put(185.7,77)   {\sse$ U_j $}
  \put(210,104)    {$ = $}
  \put(281.2,29.9) {\sse$ \ol\alpha $}
  \put(271.5,107)  {\sse$ U_i $}
  \put(305.2,168)  {\sse$ \Xd_\ka^{} $}
  \put(329.2,185.3){\sse$ \beta $}
  \put(355.5,127)  {\sse$ U_j $}
  }
Here the second equality is presented to facilitate comparison with formulas in 
\cite{ffrs5}; it follows by repeated use of functoriality of the braiding,
the various properties of $A$ (symmetric, special, Frobenius) and $X$
($A$-bimodule) and of $h^{(ij)}_\beta,h^{(ij)}_{\ol\alpha}$ (bimodule 
morphisms), and by invoking the isomorphism $\End(\one)\Cong\EndAA(A)$.

  
\section{Quantum field theory}

A main motivation for the investigations in \cite{ostr} came from work in
low-dimensional   
quantum field theory, and that area also constitutes an important arena
for applications of results like Theorem O.
This section describes some of the connections with quantum field theory.

\subsubsection*{Subfactors}

As already mentioned, braided induction first appeared in the study of 
subfactors, or more precisely, of a new construction of subfactors
that was inspired by examples from quantum field theory, see e.g.\ 
\cite{lore,xu3,izum9}. In that context it is defined in terms of the
so-called statistics operators (describing the braiding) and of the
canonical and dual canonical endomorphisms of the subfactor.

In a series of papers starting with \cite{boev,boek}, various aspects of 
braided induction were studied in this setting
(for reviews see e.g.\ \cite{boev5} or section 2 of \cite{evpi}).
The results include in particular the subfactor version of Theorem O, 
describing the fusion rules of $N'$-$N'$-morphisms for a finite index subfactor
$N\,{\subseteq}\,N'$, which was established as Theorem 6.8 of \cite{boek}.
Indeed, the formulation of Theorem O in \cite{ostr} was directly motivated by 
this result, as were several other conjectures, all of which have meanwhile 
been proven in a purely categorical setting as well.  

The proof of Theorem 6.8 in \cite{boek} is based on techniques and results 
from the theory of subfactors. As a consequence, the statement of the theorem 
itself is somewhat weaker than the one of Theorem O. More specifically, in the 
subfactor context the modular tensor category \C\ arises as a category 
$\mathcal E\hspace*{-.9pt}nd(N)$ of endomorphisms of a type III factor $N$
(see e.g.\ \cite[Def./Prop.\,2.5]{muge13}).      
This in turn implies that \C\ as well as 
$\CM\eq\mathcal H\hspace*{-.9pt}om(N,N')$ and 
$\Cs\eq\mathcal E\hspace*{-.9pt}nd(N')$ are *-categories. 
A *-category \C\ is a \complex-linear category endowed with a family
$\{ *_{U,V}{:}~\Hom(U,V)\To\Hom(V,U) \,|\, U,V\iN\obj(\C) \}$ of maps which are
antilinear, involutive and contravariant, as well as, in case \C\ is a monoidal
category, monoidal. In a *-category the dimension of any object is positive. 
This is the reason why positivity of the dimension was included in the original
form \cite{ostr} of Theorem O; as already stressed, this condition is not 
essential.

Let us also point out that in the subfactor setting algebras arise in the guise
of so-called $Q$-systems as introduced in \cite{long6}. Indeed \cite{evpi},
a $Q$-system is the same as a symmetric special Frobenius *-algebra, where the 
*-property means that $\Delta \eq{*_{\!A,A}}(m)$ and $\eps\eq{*_{\one,A}}(\eta)$.
The $Q$-system associated to $N\,{\subseteq}\,N'$ is constructed from the
embedding morphism $N\,{\hookrightarrow}\,N'$ and its two-sided adjoint.
Given a separable type III factor $N$ there is a bijection between finite index
subfactors $N\,{\subseteq}\,N'$ and $Q$-systems in 
$\mathcal E\hspace*{-.9pt}nd(N)$ (see e.g.\ \cite{muge8}).


\subsubsection*{Conformal field theory}

One of the basic ideas in the subfactor studies was to introduce
the matrix \erf{zij} and to establish that it enjoys various remarkable
properties, like commuting with certain matrices $S$ and $T$ that generate a 
representation of the modular group SL$(2,\zet)$. These properties are, in fact,
familiar from the matrix describing the so-called torus partition function in 
two-dimensional rational conformal field theory (RCFT), see e.g.\ \cite{DIms}.
One of the basic ingredients of RCFT is the `category of chiral sectors',
which provides basic input data for the field content of a given model of RCFT.
This category \C\ is a modular tensor category. In the formalization of
RCFT through conformal nets of subfactors \C\ is given by the endomorphism
category $\mathcal E\hspace*{-.7pt}nd(N)$ of the relevant factor;
in another formalization it arises as the representation category of a
conformal vertex algebra \cite{lepo12}.

However, the properties of the matrix $z$ verified in the subfactor studies are
by far not the only requirements that a valid torus partition function of RCFT 
must fulfill. (In fact, one knows of many examples of non-negative integral 
matrices in the commutant of the action of SL$(2,\zet)$ that do not describe a 
valid partition function, see e.g.\ section 4 of \cite{fusS} and \cite{gann17}.)
It has been shown that a necessary \cite{fjfrs2} and sufficent \cite{fuRs4,fjfrs}
condition for {\em all\/} constraints on the torus partition function, as well 
as on all other correlation functions of an RCFT, to be satisfied is that there 
exists a simple symmetric special Frobenius algebra in \C, through the \rep\
theory of which the field content of the theory can be understood.
This algebra is determined only up to Morita equivalence. A Morita invariant 
formulation leads precisely to the definition \erf{zij} for $z$.
  
In the approach to RCFT via Frobenius algebras $A$ in modular tensor categories
(for a brief survey see e.g.\ \cite{fuRs11,scfr2}), the role of the bimodule
category \CAA\ is to specify the allowed types of topological defect lines.
More generally, for any pair $A$, $B$ of symmetric special Frobenius algebras
the defect lines described by \CAB\ constitute one-dimensional phase boundaries 
on the two-dimensional world sheet, and their properties can be used to extract
symmetries and order-disorder dualities of the RCFT.
(The phase boundaries are referred to as {\em topological\/} defect lines
because, as it turns out \cite{ffrs5}, correlation functions remain invariant
under continuous variation of their location.)
The results from \cite{ffrs5} that were referred to and used above 
arose from the study of such defect lines in RCFT.

Fusion rules for topological defect lines in RCFT were first studied in
\cite{pezu5,pezu6},
where it was e.g.\ verified that the formula \erf{Ns=uuu} yields non-negative
integers when one inserts for the entries of the basis transformation matrix 
$u$ explicit numerical values that can be extracted from known data for 
the $\mathfrak{sl}(2)$ Wess-Zumino-Witten models,
a particular class of RCFT models. 

\medskip

The objects of the category \CA\ of $A$-modules correspond to conformal
boundary conditions of the RCFT. That \CA\ is a module category over \CAA\
corresponds in the RCFT context to the fact that one can fuse a topological
defect line with a boundary condition, thereby obtaining another boundary
conditions of the theory \cite{ffrs5}.

It is also worth mentioning that the matrix $z(A)$ naturally arises not 
only in the discussion of the bimodule category \CAA, but also of the module 
category \CA. In particular, similarly as $\dimc(\Fs)\eq {\rm tr} \big( 
z(A)^{\rm t}_{}z(A) \big)$, one finds that the number of isomorphism classes
of simple $A$-modules equals ${\rm tr} \big(z(A)\big)$
\cite[Theorem\,5.18]{fuRs4}. A crucial ingredient of the proof of this relation 
in \cite{fuRs4}\,%
  \footnote{~In the arguments in \cite{fuRs4} the three-dimensional topological
  field theory associated \cite{TUra,KAss} to $\C$ is used. But they
  can in fact be formulated entirely in categorical language, a crucial
  ingredient being semisimplicity of \C.}
is the matrix $S^A$ that implements a modular S-transformation on the conformal 
one-point blocks on the torus.  The rows of this matrix are naturally labeled 
by the isomorphism classes of simple $A$-modules, while its columns correspond 
to a basis of the subspace of so-called {\em local\/} morphisms 
\cite[Def.\,5.5]{fuRs4} in $\bigoplus_{i\in I}\Hom(A\oti U_i,U_i)$. The 
dimension of the latter space is easily seen to be ${\rm tr} \big(z(A)\big)$, 
and the result then follows by observing \cite[Eq.\,(5.113)]{fuRs4}
that the matrix $S^A$ is invertible and hence in particular square. 
(The notion of local morphisms is also closely related to certain interesting
endofunctors of \C, compare section 3.1 of \cite{ffrs}.)

In view of this connection with \CA\ it not surprising that the structure 
constants $\Ns \ka\kap\kapp$ of the bimodule fusion rules can be related to
the matrix $S^A$ as well. To see that this is the case, one can start from
the fact that \erf{pic_navf_5} is the idempotent that corresponds to forming 
the tensor product over $A$, which allows one to derive that
  \eqpic{pic_navf_9} {290} {63} {
           \put(0,-7){   
  \put(-10,5)     {\Includeourbeautifulpicture 9}
  \put(176,75)    {$ =~\dsty\sum_{\kapp\in\Is}\sum_{i\in I}\Ns\ka\kap\kapp
                       n_\kapp^i\,s_{i,j} $}
  \put(32.8,64)   {\sse$ \Xd_\ka^{} $}
  \put(55.5,75)   {\sse$ A $}
  \put(73.4,64)   {\sse$ \Xd_\kap^{} $}
  \put(96.3,97)   {\sse$ U_j $}
          }
  }
with $n_\ka^i\eq \dimc\,\Hom(U_i,\Xd_\ka)$. By a calculation repeating the steps
in formula (5.131) of \cite{fuRs4}, one then arrives for any $j\iN I$ at the sum
rule
  \be
  \sum_{\kapp\in\Is} n^\kapp_j \Ns\ka\kap\kapp
  = \sum_{i\in I} \sum_{\gamma} \sum_{\mu,\nu\in J} m_{\ka^\vee}^\mu \,
  m_{\kap}^\nu \, S^A_{\mu,i\gamma}\, \frac
  {s^{}_{i,\overline\jmath}} {s^{}_{i,0}}\, \big(S^A\big)^{-1}_{i\gamma,\nu} \,,
  \ee
where the $\gamma$-summation extends over a basis of the local morphisms in 
$\Hom(A\oti U_i,U_i)$ and the $\mu$- and $\nu$-summations over a set $J$ of
representatives of simple left $A$-modules, while the integers 
$m_\ka^\mu \eq \dimc\,\HomA(M_\mu,X_\ka)$ count the multiplicity of the simple 
left $A$-module $M_\mu$ in the simple $A$-bimodule $X_\ka$,
regarded as a (not necessarily simple) left $A$-module.

\vskip2em

\noindent {\small {\bf Acknowledgements} \\[1pt]
JF is partially supported by VR under project no.\ 621--2003--2385, IR by the 
EPSRC First Grant EP/E005047/1 and the PPARC rolling grant PP/C507145/1, and
CS by the Collaborative Research Centre 676 ``Particles, Strings and the
Early Universe - the Structure of Matter and Space-Time''.}

 \vskip 3em

 \newcommand\wb{\,\linebreak[0]} \def\wB {$\,$\wb}
 \newcommand\Bi[2]    {\bibitem[#2]{#1}}
 \newcommand\BOOK[4]  {{\sl #1\/} ({#2}, {#3} {#4})}
 \newcommand\JO[6]    {{\em #6}, {#1} {#2} ({#3}), {#4--#5} }
 \newcommand\inBO[9]  {{\em #8}, in:\ {\sl #1}, {#2}\ ({#3}, {#4} {#5}),
                      p.\ {#6--#7} {{\tt [#9]}}}
 \newcommand\J[7]     {{\em #7}, {#1} {#2} ({#3}), {#4--#5} {{\tt [#6]}}}
 \newcommand\Pret[2]  {{\em #2}, pre\-print {\tt #1}}
 \def\jf    {J.\ Fuchs}
 \def\adma  {Adv.\wb Math.}
 \def\anma  {Ann.\wb Math.}
 \def\coma  {Con\-temp.\wb Math.}
 \def\comp  {Com\-mun.\wb Math.\wb Phys.}
 \def\jpaa  {J.\wB Pure\wB Appl.\wb Alg.}
 \def\nupb  {Nucl.\wb Phys.\ B}
 \def\phlb  {Phys.\wb Lett.\ B}
 \def\pnas  {Proc.\wb Natl.\wb Acad.\wb Sci.\wb USA} 
 \def\rvmp  {Rev.\wb Math.\wb Phys.}
 \def\trgr  {Trans\-form.\wB Groups}
 \def\taac  {Theo\-ry\wB and\wB Appl.\wB Cat.}
 \def\CUP    {{Cambridge University Press}}
 \def\EMS    {{European Mathematical Society}}
 \def\SV     {{Sprin\-ger Ver\-lag}}
 \def\Ca     {{Cambridge}}
 \def\NY     {{New York}}

\end{document}